\documentclass[11pt]{article}
\usepackage[centertags]{amsmath}
\usepackage{amssymb,amsfonts,theorem,stmaryrd}
\usepackage[all]{xy}

\allowdisplaybreaks

\theorembodyfont{\slshape}
\newtheorem{thm}{Theorem}[section]
\newtheorem{prop}[thm]{Proposition}

\newtheorem{lemma}[thm]{Lemma}
\newtheorem{cor}[thm]{Corollary}
\theorembodyfont{\rmfamily}

\newenvironment{pf}{\textit{Proof.}}{\hfill$\boxbox$}

\newcommand{\R}{\mathbb{R}}
\newcommand{\N}{\mathbb{N}}
\newcommand{\Z}{\mathbb{Z}}
\newcommand{\C}{\mathbb{C}}
\renewcommand{\H}{\mathbb{H}}

\newcommand{\spc}{\mathcal{N}}

\newcommand{\dom}{\operatorname{dom}}
\newcommand{\spec}{\operatorname{spec}}

\renewcommand{\min}{min}
\renewcommand{\max}{max}

\newcommand{\SL}{\operatorname{SL}}

\newcommand{\dev}{\operatorname{dev}}
\newcommand{\hol}{\operatorname{hol}}

\newcommand{\supp}{\operatorname{supp}}
\newcommand{\Defo}{\operatorname{Def}}
\newcommand{\grad}{\operatorname{grad}}
\newfont{\hugemath}{cmsy10 scaled 3000}

\begin{document}
\title{The Laplacian on hyperbolic 3-manifolds with Dehn surgery type singularities}
\author
{Frank Pf\"aff\-le \\ Universit\"at Potsdam \and Hartmut Wei{\ss} \\ LMU M\"unchen}

\maketitle

\begin{abstract}
We study the spectrum of the Laplacian on hyperbolic 3-manifolds
 with Dehn surgery type singularities and its dependence on the generalized
 Dehn surgery coefficients. 
\end{abstract}

\tableofcontents

\section{Introduction}
Let $M$ be a complete non-compact hyperbolic 3-manifold of finite volume. 
The Laplacian considered as a symmetric densely-defined operator  
\[
\Delta: C_0^{\infty}(M) \rightarrow L^2(M)
\]
is essentially selfadjoint, cf.~\cite{Gaf}. The essential spectrum of the
unique selfadjoint extension of $\Delta$ consists of the interval
$[1,\infty)$, cf.~\cite{DL}, \cite{MP}. Let us further denote by $0=\lambda_0 <\lambda_1 \leq \ldots
\leq \lambda_k < 1$ the finitely many eigenvalues below the essential
spectrum, cf.~\cite{LP}.  

If $(M_i)_{i \in \N}$ is a sequence of compact hyperbolic $3$-manifolds, which
converges to $M$, say in the pointed Lipschitz topology, one can ask in what
sense the spectrum of the limit manifold is related to the spectra of the
approximators. Note that since $M_i$ is a compact manifold, the spectrum of
the Laplacian on $M_i$ is discrete for each $i \in \N$. Let us assume for
simplicity that $M$ has a single rank-2 cusp. 
The following results are known:

In \cite{CC}, B.~Colbois and G.~Courtois show that the eigenvalues of $M$
below the essential spectrum are limits of eigenvalues of the $M_i$. 
More precisely, if $0 = \lambda_0^i < \lambda_1^i \leq \ldots \leq \lambda_{k(i)}^i <1 $ are
the 
eigenvalues of $M_i$ smaller than $1$, then for $i$ large enough one has
$k(i)\geq k$ and further
\begin{equation}\label{Convergence}
\lim_{i \rightarrow \infty} \lambda_j^i = \lambda_j
\end{equation}
for $j=0, \ldots, k$.

In \cite{CD}, I.~Chavel and J.~Dodziuk show that the eigenvalues of the $M_i$
accumulate in the interval $[1,\infty)$ as $i \rightarrow \infty$. 
Moreover, they determine the precise rate of clustering in terms of geometric
data of the degenerating tube. 
Namely, with $\spc_{\Delta_i,M_i}[1,1+x^2]=\vert\{\lambda \in \spec \Delta_i :
1 \leq \lambda \leq 1+x^2\}\vert$ 
denoting the spectral counting function of the Laplacian on $M_i$, they obtain
the estimate 
\begin{equation}\label{Clustering}
\spc_{\Delta_i,M_i}[1,1+x^2] = \frac{x}{2 \pi} \log \left(\frac{1}{l_i}\right)
+ O_x(1)\,, 
\end{equation}
where $l_i$ is the length of the shortest closed geodesic in $M_i$.

Let us mention that similar questions have been studied for the Laplacian on
differential forms by J.~Dodziuk and J.~McGowan in \cite{DMc} and for the
Dirac operator on complex spinors by C.~B\"ar in \cite{Bae} and F.~Pf\"aff\-le
in \cite{Pf}.  

From the point of view of the deformation theory of $M$ it is natural to bring
a certain class of singular hyperbolic 3-manifolds into play: While due to
Mostow-Prasad rigidity, cf.~\cite{Mos}, \cite{Pra}, the hyperbolic structure
on $M$ may not be deformed through complete hyperbolic structures, there is
actually a real 2-dimensional deformation space of incomplete hyperbolic
structures parametrized by the so-called generalized Dehn surgery
coefficients. This is the essence of Thurston's Hyperbolic Dehn Surgery
Theorem, cf.~\cite{Th}. The existence of sequences $M_i$ as above is in fact a
consequence of that theorem. 

These incomplete structures are said to have Dehn surgery type singularities,
special cases include hyperbolic cone-manifold structures and in particular
smooth hyperbolic structures on certain topological fillings.  

The aim of this article is to study basic spectral properties of the
Laplacian on hyperbolic 3-manifolds with Dehn surgery type singularities and
to prove analogues of the asymptotic statements in equations (\ref{Convergence})
and (\ref{Clustering}) for this wider class of hyperbolic manifolds. 
The main results are Theorems~\ref{convergence} and \ref{spec_count}. Note
that due to Corollary \ref{extends_results} our results include the results of
\cite{CC} and \cite{CD}. Moreover, we emphasize a "continuous" aspect of these phenomena, namely all our estimates take place on the deformation space of structures (and do not make reference to specific sequences of manifolds).

The authors would like to thank Christian B\"ar for useful conversations and
SFB 647 "Raum - Zeit - Materie" for financial support, furthermore Prof.~Kalf for pointing out the references \cite{Sea1} and \cite{Sea2} to us.

\section{Hyperbolic Dehn surgery}\label{Dehn_surgery}
Let $M$ be a complete hyperbolic 3-manifold of finite volume. As a consequence
of the Margulis Lemma, cf.~\cite{KM}, $M$ has only finitely many ends all of
which are rank-2 cusps. More precisely, if for $\mu >0$ we look at the thick-thin decomposition of $M$, i.e.~$M = M_{(0,\mu)} \cup
M_{[\mu,\infty)}$, where  
$$
M_{(0,\mu)}=\{ p \in M : inj_p < \mu \}
$$
is the $\mu$-{\em thin} part of $M$, and
$$
 M_{[\mu,\infty)} =\{ p \in M : inj_p \geq \mu\} 
$$
the $\mu$-{\em thick} part, then the Margulis Lemma asserts the existence of a universal constant $\mu_0$, such that the components of
the $\mu$-thin part of $M$ have standard geometry for $\mu<\mu_0$: They are either {\em rank-2
  cusps} or {\em smooth tubes}.  

A {\em rank-2 cusp} is the quotient of a horoball in $\H^3$ by a rank-2 free
abelian group of parabolic isometries, which we in the following denote by
$\Gamma_{cusp}$. The intrinsic geometry of the boundary horosphere is that of
flat $\R^2$, its principal curvatures are constantly $1$. Let
$T^2_{cusp}=\R^2/\Gamma_{cusp}$ be the corresponding flat torus. Then the cusp
based on $T^2_{cusp}$ is given as the Riemannian manifold 
$$
((0,\infty) \times T^2,dt^2 + e^{-2t}g_{T^2})\,,
$$
where $t \in (0,\infty)$ and $g_{T^2}$ denotes the flat metric on $T^2$.

A {\em smooth tube} is the quotient of the distance tube of a geodesic
$\gamma$ in $\H^3$ by an infinite cyclic group of hyperbolic isometries. Let
us denote its generator in the following by $\phi$. The isometry $\phi$ is a
screw-motion along $\gamma$, i.e.~after orienting $\gamma$, we can associate
the translation length $l > 0$ and the rotation angle $t \in [0,2\pi)$ with
  $\phi$. 
The boundary of the distance tube of radius $r$ is intrinsically flat, its
principal curvatures are given by $\coth(r)$ and $\tanh(r)$. 
Let  $\Gamma$ be the lattice in $\R^2$ generated by the vectors $(2\pi,0)$ and
$(t,l)$ and let $T^2=\R^2/\Gamma$. 
Then the tube of radius $R$ based on  $T^2$ is given as the Riemannian manifold 
$$
((0,R) \times T^2, dr^2 +\sinh(r)^2 d\theta^2 + \cosh(r)^2dz^2)\,,
$$
where r $\in (0,R)$ and $(\theta,z) \in \R^2$.

Let $\bar{M}$ denote the compact core of $M$ obtained by removing the cusp
components from the $\mu$-thin part of $M$. 
The boundary components of $\bar{M}$ are horospherical tori, whose injectivity
radii satisfy a universal lower bound. 
For simplicity we will assume in the following that $M$ has a single cusp.   

Let $\Defo(M)$ be the deformation space of (possibly incomplete) hyperbolic
structures on $M$, for precise definitions see \cite{CHK}. 
A hyperbolic structure is determined by a local diffeomorphism $\dev:
\tilde{M} \rightarrow \H^3$, the {\em developing map}, which is equivariant
w.r.t.~a group homomorphism $\hol: \pi_1M \rightarrow \SL_2(\C)$, the {\em
  holonomy representation}. 
The topology of $C^{\infty}$-convergence on compact subsets of $\tilde{M}$ on
the space of developing maps induces a topology on $\Defo(M)$. 
The map obtained by sending a hyperbolic structure to its holonomy
representation induces a local homeomorphism $\Defo(M) \rightarrow
X(\pi_1M,\SL_2(\C))$, where the latter is the space of group homomorphisms
$\rho:\pi_1 M \rightarrow \SL_2(\C)$ considered up to conjugation in $\SL_2(\C)$, cf.~\cite{Gol}, see also \cite{CHK}. 

We fix generators $\mu, \lambda \in \pi_1 \partial \bar{M}$. Let $\rho_0$ be
the holonomy of the complete structure and $\chi_0$ the corresponding element in $X(\pi_1M,\SL_2(\C))$. 
For $\rho$ a deformation of $\rho_0$ we consider the complex lengths
$\mathcal{L}_{\mu}(\rho)$ and $\mathcal{L}_{\lambda}(\rho)$ of the isometries
corresponding to $\rho(\mu)$ and $\rho(\lambda)$. Since the complex length is
invariant under conjugation, but only determined up to addition of multiples
of $2\pi i$ and multiplication by $\pm1$, $\mathcal{L}_{\mu}$ and
$\mathcal{L}_{\lambda}$ may be considered as multi-valued functions on
$X(\pi_1 \bar{M}, \SL_2(\C))$.

In \cite{BP} it is shown that on a branched cover of a neighbourhood of
$\chi_0$ the functions $\mathcal{L}_{\mu}$ and $\mathcal{L}_{\lambda}$ can be
defined as single-valued functions. 
More precisely, there exist open sets $0\in U \subset \C$ and $\chi_0 \in V
\subset X(\pi_1 \bar{M}, \SL_2(\C))$, a branched cover $\pi: U \rightarrow V$
and a holomorphic map $f: U \rightarrow \C$ such that 
$$
\mathcal{L}_{\mu}(\pi(z)) = z
$$ 
and
$$
\mathcal{L}_{\lambda}(\pi(z))=f(z) \,.
$$
Furthermore for $z$ small enough (w.l.o.g.~for $z \in U$), the equation
\begin{equation}\label{genDehnsurgcoeff}
x z + yf(z)=2\pi i 
\end{equation}
has a unique solution $(x,y) \in \R^2 \cup \{ \infty \}$, the so-called {\em
  generalized Dehn surgery coefficients}. 
The map $f$ satisfies $f(-u)=-f(u)$ (in particular $f(0)=0$), such that the
  generalized Dehn surgery coefficients are determined by a character $\chi
  \in V$ up to sign. 
Hence for a sufficiently small neighbourhood $W$ of the complete structure in
  $\Defo(M)$ there is a well-defined map
$$
DS: W \rightarrow (\R^2 \cup \{\infty\}) / \pm 1\,.
$$
With these preparations we can state Thurston's Hyperbolic Dehn Surgery
Theorem, cf.~\cite{Th}, see also \cite{BP}, \cite{CHK}: 
\begin{thm}[W.P.~Thurston]\label{dehn_surgery} There is a neighbourhood $W
  \subset \Defo(M)$ of the complete structure such that the map 
$$
DS: W \rightarrow (\R^2 \cup \{\infty\})/\pm 1
$$
is a homeomorphism onto a neighbourhood $W'$ of $\infty \in (\R^2 \cup
\{\infty\})/\pm 1$. 
\end{thm}
In the following we give a geometric description of how these incomplete
structures look like:  

Let $g$ denote the complete hyperbolic metric on $M$. Recall that $\partial
\bar{M}$ equipped with $g$ is a horospherical torus, i.e.~intrinsically flat
with principal curvatures $1$. 
A structure which is close to the complete structure may be represented by a
hyperbolic metric $g'$ on $\bar{M}$, which is $C^{\infty}$-close to $g$
restricted to $\bar{M}$, such that $\partial \bar{M}$ equipped with $g'$
becomes intrinsically flat with principal curvatures $\coth(R)$ and $\tanh(R)$
for $R>0$ large. 
Then a {\em singular tube} of radius $R$ is added, such that principal
curvature lines match. 

A {\em singular tube} of radius $R$ is the following obvious generalization of
a smooth tube as above: Let $\Gamma \subset \R^2$ be any lattice and consider
the metric 
\begin{equation}\label{singtubemetric}
dr^2 + \sinh(r)^2d\theta^2 + \cosh(r)^2dz^2
\end{equation}
on $(0,R) \times \R^2/\Gamma$, where $r \in (0,R)$ and $(\theta,z) \in \R^2$. 
We denote the tube of radius $R$ based on the torus $T^2=\R^2/\Gamma$ by
$T^2_{(0,R)}$. 
The following cases occur:
\begin{enumerate}
\item $\Gamma \cap \{z=0\} \neq \{0\}$: 
In this case $\Gamma$ is spanned by unique vectors $(\alpha,0)$ with
$\alpha>0$ and $(t,l)$ with $0\leq t<\alpha$ and $l>0$. Then
$T^2_{(0,R)}$ is the smooth part of a cone tube with cone angle
$\alpha$, length $l$ and twist $t$. 
In the special case $\alpha = 2\pi$ this is nothing but a smooth tube with the
core geodesic removed. 
The principal curvature lines corresponding to $\coth(R)$ close up and are
isotopic to a curve $p \mu + q \lambda$ with $\mu,\lambda \in
\pi_1\partial\bar{M}$ as above and $p,q$ coprime integers. 
The generalized Dehn surgery coefficients of the structure are given by
$(x,y)=\frac{2\pi}{\alpha}(p,q)$.     
\item $\Gamma \cap \{z=0\} = \{0\}$: In this case the metric completion is not
  a manifold. 
The generalized Dehn surgery coefficients of the structure are of the form
  $(x,y)$ with $x/y$ irrational. 
\end{enumerate}
We will generally say that a hyperbolic 3-manifold obtained in this way has
Dehn surgery type singularities.

The precise shape of the deformed tube is determined by the deformed holonomy
of the boundary torus, i.e.~by the complex lengths $\mathcal{L}_{\mu}$ and
$\mathcal{L}_{\lambda}$. 
These in turn are determined by the generalized Dehn surgery coefficients via
equation \ref{genDehnsurgcoeff}, however this dependence is not explicit.
 Qualitatively we can say the following:
\begin{lemma}\label{tube_geom}
Let $M^{\infty}$ be a complete hyperbolic 3-manifold of finite volume with a single
cusp. 
Then:
\begin{enumerate}
\item For any $\varepsilon >0$ there is a neighbourhood $W \subset \Defo(M^{\infty})$
  of the complete structure such that $diam\, T^2 < \varepsilon$ for any
  hyperbolic structure in $W$, where $T^2=\R^2/\Gamma$ is the base of the
  corresponding tube.   
\item For any neighbourhood $W \subset \Defo (M^{\infty})$ of the complete structure
  there are constants $C_1,C_2 >0$ such that $C_1 \leq e^{2R} area\, T^2 \leq
  C_2$ for any hyperbolic structure in $W$. 
\end{enumerate}
\end{lemma}
\begin{pf}
Recall that the complex lengths $\mathcal{L}_{\mu}$ and
$\mathcal{L}_{\lambda}$ are defined as single valued functions $z$ and $f(z)$
on a branched cover of a neighbourhood of the complete structure with $f$
holomorphic and $f(0)=0$. 
Hence $\Gamma$ is spanned by arbitrarily short
vectors for hyperbolic structures close enough to the complete one. 

Let $T^2_r$ denote the torus $T^2=\R^2/\Gamma$ equipped with the Riemannian
metric $\sinh(r)^2d\theta^2+\cosh(r)^2dz^2$. 
Then $area\,T^2_r =\sinh(r)\cosh(r) area \,T^2$ and $area \, T^2_R$ differs
from the area of the horospherical torus $\partial\bar{M}^{\infty}$ by some bounded
amount depending on the neighbourhood $W \subset \Defo(M^\infty)$. 
Since $\sinh(r)\cosh(r) \sim e^{2r}$, the second claim follows. 
\end{pf}\\
\\
Note that a cone tube based on $T^2$ with small diameter (and hence small
area) may have small cone angle or not. 
In any case it is easy to see that the length of the tube has to be small. 

Further, for a cone tube one has $area\, T^2 = l\alpha$, where $\alpha$ is the
cone angle and $l$ the length of the tube. 
In particular, if we restrict to smooth fillings, i.e.~$\alpha=2\pi$, Lemma
\ref{tube_geom} gives us constants $C_1',C_2'>0$ such that 
$$
C_1' \leq e^{2R} l \leq C_2' \Leftrightarrow R-C_1'' \leq \frac{1}{2} \log
\frac{1}{l} \leq R + C_2'' \,,  
$$ 
which is the estimate used in \cite{CD} and \cite{Bae}. 
This suggests that $area \,T^2$ in our arguments should play the role of $l$
in the arguments of \cite{CD} and \cite{Bae}.

\begin{lemma}\label{lem_minftydef}
Consider $c> 4$ and $\beta>0$.
Then there exists $\mu>0$ below the Margulis constant $\mu_0$
and a neighbourhood $W$ of $M^\infty\in\Defo(M^\infty)$ such that
one has:
\begin{enumerate}
\item For any $M\in W$ the $\mu$-thick part
  $M_{[\mu,\infty)}$ is $(1+\beta)$-quasi-isometric to
  $M^\infty_{[\mu,\infty)}$, 
\item any $M\in W$ contains a singular tube $T^2_{(0,\rho+c]}$ for some $\rho >0$
  such that $T^2_{[\rho,\rho+c]}\subset M_{[\mu,\infty)}$.
\end{enumerate}
\end{lemma}
\begin{pf}
The proof is evident from the discussion above.
\end{pf}

\section{The spectrum of the Laplacian}
Let $M$ be a hyperbolic 3-manifold with Dehn surgery type singularities. 
The Laplacian $\Delta$ considered on $\dom \Delta=C_0^{\infty}(M)$ is a
symmetric, densely defined operator in $L^2(M)$. 
Since $M$ is incomplete, 
we cannot expect the Laplacian to be essentially selfadjoint on that domain;
in fact we will see that it never is. 
Nevertheless, since $\Delta$ is nonnegative on $C_0^{\infty}(M)$,
i.e.~$\langle \Delta f,f\rangle_{L^2} \geq 0$ for all $f\in C_0^{\infty}(M)$,
there is always a distinguished selfadjoint extension at hand, namely the
so-called Friedrichs extension of $\Delta$. 

In the following we briefly review the construction of the Friedrichs
extension of a semibounded symmetric operator. 
Recall that if $\mathcal{H}$ is a Hilbert space and $q$ a quadratic form
defined on a dense domain in $\mathcal{H}$ such that $q\geq c$ for some $c \in
\R$ and $q$ is closable with closure $\bar{q}$, then there exists a unique
selfadjoint operator $A$ with $\dom A \subset \dom \bar{q}$ and $\langle A f,g
\rangle = q(f,g)$ for all $f \in \dom q \cap \dom A$ and $g \in \dom
q$. Furthermore, $A$ satisfies the same lower bound as $q$, i.e.~$A\geq c$. 
The domain of $A$ is given by 
$$
\dom A = \{ f \in \dom \bar{q} : \exists \,h \in \mathcal{H} \text{ s.t. }
\bar{q}(f,g)=\langle h,g \rangle \,\forall g \in \dom q \} 
$$
and then $Af=h$.

The spectral theorem for selfadjoint operators yields the well-known
variational characterization of the eigenvalues of $A$ below the essential
spectrum:
\begin{thm}\label{varcharacterization}
For 
$
\lambda_k=\inf\limits_{\substack{ V\subset \dom q \\ \dim(V)=k}}\;
 \sup\limits_{ f\in V\setminus\{ 0\}}\; \frac{q(f,f)}{|f|^2}
$
one has:
\begin{enumerate}
\item the sequence $(\lambda_k)$ is non-decreasing, and
  $\lambda_k\to\lambda_\infty\le\infty$, 
\item the $\lambda_n<\lambda_\infty$ are precisely the eigenvalues of $A$ below
  $\lambda_\infty$, 
\item $\lambda_\infty$ is the bottom of the essential spectrum of $A$.
\end{enumerate}
\end{thm}

Now if $A_0$ is a densely defined symmetric operator in $\mathcal{H}$ with
$A_0 \geq c$, then $q(f,g)=\langle A_0f,g\rangle$ for $f,g \in \dom A_0$ is
closable and trivially $q \geq c$. One has $\dom A = \dom A_0^* \cap \dom
\bar{q}$ and $Af=A_0^*f$ for the corresponding selfadjoint operator $A$, which
is the so-called {\em Friedrichs extension} of $A_0$. 

If $M$ is a Riemannian manifold (without boundary) we may apply this
construction to the Dirichlet energy $q(f) = \int_M |df|^2$ on $\dom q =
C_0^{\infty}(M)$ to obtain the Friedrichs extension of $\Delta$, in the
following denoted by $\Delta_{Fr}$. 
One has $\dom\bar{q}=\dom d_{min}$ and $\dom \Delta_{Fr}=\dom \Delta_{max}
\cap \dom d_{min}$, where for a differential operator $P$ acting on compactly
supported smooth sections of some vector bundle (equipped with a Euclidean
metric) we set 
$$
P_{max} = (P^t)^* \quad \text{ and } \quad P_{min} = \bar{P}\,.
$$
Here $P^t$ denotes the formal adjoint of $P$ and $\bar{P}=P^{**}$ the closure
of $P$. 
It is easy to see that
$$
\dom P_{\max} = \{ s \in L^2 : Ps \in L^2 \} 
$$
and 
$$
\dom P_{\min} = \{ s \in L^2: \exists s_n \in C_0^{\infty} \text{ with } s_n
\overset{L^2}{\rightarrow} s,\, s_n \overset{L^2}{\rightarrow} Ps \} \,,
$$   
where $P$ is applied to $L^2$-sections in a distributional sense.

For example, if $M$ is the interior of a compact manifold with boundary
$\hat{M}$, i.e.~$M=\hat{M} \setminus \partial\hat{M}$, then $\dom
\bar{q}=H^1_0(M)$ and the Friedrichs extension of $\Delta$ is the Dirichlet
extension with  
$$
\dom \Delta_{Dir} = H^2(M) \cap H^1_0(M)\,.
$$
If we apply the Friedrichs construction to $q(f)=\int_M |df|^2$ on $\dom q =
C^{\infty}(\hat{M})$, then $\dom \bar{q} = H^1(M)$ and the corresponding
selfadjoint operator is the Neumann extension of $\Delta$ with 
$$
\dom \Delta_{Neu} = \{ f \in H^2(M) : \nu f \in H^1_0(M) \}\,,
$$
where $\nu$ is a smooth extension of a unit normal to the boundary.

Returning to $M$ being a hyperbolic 3-manifold with Dehn surgery type
singularities, let $d: C^\infty(M) \rightarrow \Omega^1(M)$ denote the
exterior differential on functions and $\delta: \Omega^1(M) \rightarrow
C^\infty(M)$ the divergence on $1$-forms. As usual let $H^1(M)=\dom d_{max}$
and $H^1_0(M)=\dom d_{min}$. From the preceding discussion it is clear that
$\Delta_{Fr} = \delta_{max}d_{min}$. 
\begin{thm}[$L^2$-Stokes]\label{L2_stokes} Let $M$ be a hyperbolic 3-manifold
  with Dehn surgery type singularities. 
Then for $f \in \dom d_{max}$ and $w \in \dom \delta_{max}$ one has
$$
\langle df, \omega \rangle_{L^2} = \langle f, \delta \omega \rangle_{L^2}\,.
$$
\end{thm}
\begin{pf}
We may w.l.o.g.~assume that $f$ and $\omega$ are smooth, i.e.~$f \in
C^{\infty}(M)$ with $\|f\|_{L^2}<\infty$ and $\|df\|_{L^2} < \infty$ and $
\omega \in \Omega^1(M)$ with $\|w\|_{L^2}<\infty$ and $\|\delta \omega\|_{L^2} <
\infty$. 
We may also replace the hyperbolic metric whose restriction to the singular
tube is the flat metric $dr^2 + r^2 d\theta^2+dz^2$, since the $L^2$-Stokes
property is unaffected by passing to a quasi-isometric metric. 
Now
$$
\int_{M\setminus 
T^2_{(0,r)}}  df \wedge \ast \omega = \int_{M\setminus T^2_{(0,r)}} f
\ast \delta \omega \pm \int_{T^2_r} f\ast \omega\,, 
$$
where $T^2_r$ denotes the cross-section $\{r\} \times T^2 \subset
T^2_{(0,R)}$. 
We wish to show that the boundary integral vanishes in the limit as $r
\rightarrow 0$ (or at least for a sequence $r_n \rightarrow 0$ as $n
\rightarrow \infty$). If we write $\omega=\phi_r dr+\phi_\theta d\theta+\phi_z
dz$, then the resriction of $\ast \omega$ to $T^2_r$ equals $\phi_r rd\theta
\wedge dz$ and  
$$
\Bigl|\int_{T^2_r}f\ast\omega\Bigr| \leq \Bigl( \int_{T^2} f^2 \, rd\theta
\wedge dz \Bigr)^{\frac{1}{2}} \cdot \Bigl( \int_{T^2} \phi_r^2 \, rd\theta
\wedge dz \Bigr)^{\frac{1}{2}}\,. 
$$
If $\varphi \in L^1(0,1)$ then Lemma 1.2 in \cite{Che} shows that there exists
a sequence $r_n \rightarrow 0$ such that $\varphi(r_n) = o(r_n^{-1}|\log
r_n|^{-1})$. 
Applied to the second factor this yields a sequence $r_n \rightarrow 0$ such
that  
$$
\Bigl( \int_{T^2} \phi_r^2 \, r_nd\theta \wedge dz \Bigr)^{\frac{1}{2}}=
o(r_n^{-\frac{1}{2}}|\log r_n|^{-\frac{1}{2}})\,. 
$$
To achieve a better estimate for the first factor we use that also
$\|df\|_{L^2}< \infty$, cf.~Lemma 2.3 in \cite{Che}: 
$$
\Bigl(\int_{T^2}\Bigl(\int_r^1 \frac{\partial f}{\partial r}(s)ds \Bigr)^2 \,
rd\theta\wedge dz \Bigr)^{\frac{1}{2}} \leq r^{\frac{1}{2}}|\log
r|^{\frac{1}{2}} \int_{T^2}\Bigl(\int_r^1 \Bigl(\frac{\partial f}{\partial
  r}(s)\Bigr)^2 s ds\Bigr) d\theta\wedge dz  
$$
Now
$$
\int_r^1 \frac{\partial f}{\partial r}(s,\theta,z)ds =
f(1,\theta,z)-f(r,\theta,z) 
$$
such that
$$
\Bigl( \int_{T^2} f^2 \, rd\theta \wedge dz
\Bigr)^{\frac{1}{2}}=O(r^{\frac{1}{2}}|\log r|^{\frac{1}{2}})
$$
Altogether we obtain that $\lim_{n \rightarrow
  \infty}\int_{T^2_{r_n}}f\ast\omega = 0$. 
\end{pf}\\
\\
This proof also shows the following Green's formula.
\begin{cor}\label{L2_green}
If we denote the singular tube in $M$ again by $T^2_{(0,r)}$, then $-\partial_r$ is the exterior normal vector field of $X$ along $T^2_r$,
and one has for 
all $f\in C^\infty(M)$ with~$\|f\|_{L^2}<\infty$,~$\|df\|_{L^2}<\infty$
and~$\|\Delta f\|_{L^2}<\infty$:
\[
\int_{T^2_{(0,r)}} |df|^2- \int_{T^2_{(0,r)}} f\cdot\Delta
f=\pm\int_{T^2_r} f\cdot\partial_r f,
\]
where all integrals are taken with respect to the induced volume measures.
\end{cor}
\begin{cor}\label{max=min}
Let $M$ be a hyperbolic 3-manifold with Dehn surgery type singularities. Then:
\begin{enumerate}
\item $d_{max}=d_{min}$, i.e.~$H^1(M)=H^1_0(M)$.
\item $0 \in \spec \Delta_{Fr}$.
\end{enumerate}
\end{cor}
\begin{pf}
Theorem \ref{L2_stokes} shows that $d_{max}=\delta_{max}^*$, where
$\delta_{max}^*$ denotes the Hilbert space adjoint of $\delta_{max}$. 
Since in general one has $\delta_{max}^*=d_{min}$, the first assertion
follows. Since $\Delta_{Fr} = \delta_{max}d_{min}$ it is enough to show that
$1 \in H^1_0(M)$. Clearly $1 \in H^1(M)$, so the second assertion follows from
the first. 
\end{pf}
\begin{cor}\label{extends_results}
Let $\hat{M}$ be a compact hyperbolic 3-manifold and let $\hat{\Delta}$ denote
the Laplacian on $\hat{M}$. 
If $M = \hat{M} \setminus \gamma$ for $\gamma \subset \hat{M}$ a closed
geodesic, then $\Delta_{Fr}$ coincides with the unique selfadjoint extension
of $\hat{\Delta}$, i.e.~$\dom \Delta_{Fr}=H^2(\hat{M})$. 
In particular, $\spec \Delta_{Fr} = \spec \hat{\Delta}$. 
\end{cor}
\begin{pf}
Clearly $H^1_0(M) \subset H^1(\hat{M}) \subset H^1(M)$, hence by Corollary
\ref{max=min} one has in particular $H^1_0(M) = H^1(\hat{M})$. 
Applying the Friedrichs construction to both form domains yields the result.
\end{pf}\\
\\
To investigate further properties of $\Delta$ we use a separation of variables
argument on the singular tube. 
A direct calculation shows that the Laplacian on $T^2_{(0,R)}$ is given by
$$
\Delta =-\partial_r^2-2\coth(2r)\partial_r +L(r)
$$
where
$$
L(r)=-\frac{1}{\sinh(r)^2}\partial_{\theta}^2-\frac{1}{\cosh(r)^2}\partial_z^2
$$
is the Laplacian on the cross-section $T^2_r$ with the induced metric. 
The volume form on $T^2_{(0,R)}$ is given by $\sinh(r)\cosh(r)dr\wedge
d\theta \wedge dz$, hence 
\begin{align*}
L^2(T^2_{(0,R)}) & \rightarrow L^2((0,R) \times T^2))\\
f & \mapsto \sinh(r)^{\frac{1}{2}}\cosh(r)^{\frac{1}{2}}f
\end{align*}
is a unitary operator. 
Since
\begin{align*}
&  \sinh(r)^{\frac{1}{2}}\cosh(r)^{\frac{1}{2}}
  \left((-\partial_r^2-2\coth(2r)\partial_r)  
\left(\sinh(r)^{-\frac{1}{2}}\cosh(r)^{-\frac{1}{2}}f\right)\right)\\ 
=& -\partial_r^2f + \left(2-\coth(2r)^2\right)f
\end{align*}
we obtain that the Laplacian on $T^2_{(0,R)}$ is unitarily equivalent to
the operator  
$$
-\partial_r^2 + 2-\coth(2r)^2+L(r)\,.
$$
If $T^2=\R^2/\Gamma$, let $\Lambda$ be the dual lattice to $\Gamma$, i.e.
$$
\Lambda = \{ \lambda \in \R^2 : \langle \lambda,\gamma \rangle \in \Z
\;\forall \, \gamma \in \Gamma \}\,. 
$$
The dual lattice $\Lambda$ is spanned by the vectors
$$
\left( \frac{l}{\sinh(r) \cdot covol(\Gamma)},  \frac{-t}{\cosh(r) \cdot
  covol(\Gamma)} \right) 
$$
and
$$
\left( \frac{-s}{\sinh(r) \cdot covol(\Gamma)}, \frac{\alpha}{\cosh(r) \cdot
  covol(\Gamma)} \right) \,. 
$$
We decompose as orthogonal Hilbert sums
$$
L^2(\R^2/\Gamma) =  \bigoplus_{\lambda \in \Lambda} \langle \Psi_{\lambda}
\rangle 
$$
and 
$$
L^2((0,R),L^2(\R^2/\Gamma)) = \bigoplus_{\lambda \in \Lambda} L^2(0,R) \otimes
\langle \Psi_{\lambda} \rangle \,, 
$$
where $\Psi_{\lambda}(x)=e^{2\pi i \langle \lambda,x \rangle} /
\sqrt{covol(\Gamma)}$ and $\langle \Psi_{\lambda} \rangle$ denotes the span of
$\Psi_{\lambda}$ in $L^2(\R^2/\Gamma)$. 
We have
$$
\spec L(r) = \left\{ (2\pi)^2 \left(\frac{\lambda_1^2}{\sinh(r)^2} + \frac
      {\lambda_2^2}{\cosh(r)^2} \right) : \lambda=(\lambda_1,\lambda_2) \in
      \Lambda \right\} 
$$
and therefore the action of the operator $-\partial_r^2 + 2-\coth(2r)^2+L(r)$
on its domain intersected with $L^2(0,R) \otimes \langle \Psi_{\lambda}
\rangle$ is given by the Schr{\"o}dinger operator 
$
P_{V_\lambda}=-\partial_r^2 + V_{\lambda}
$
with potential 
$$
V_{\lambda}(r) = 2 - \coth(2r)^2 + (2\pi)^2
\left(\frac{\lambda_1^2}{\sinh(r)^2} + \frac {\lambda_2^2}{\cosh(r)^2} \right)
\,. 
$$

Let in the following $P_V=-\partial_r^2+V(r)$ be a Schr\"odinger operator on the interval $(0,R)$ with continuous potential $V$. It has been shown by D.B.~Sears in \cite {Sea1} that if 
$$
V(r) \geq \frac{3}{4r^2}+A
$$
for some constant $A \in \R$, then the limit-point case holds at $0$ (using Weyl's classical terminology), whereas the limit-circle case holds at $0$ if
$$
|V(r)| < \frac{\delta}{r^2} +A
$$
for $0\leq \delta < 3/4$ and some constant $A \in \R$.

The following lemma asserts that it is indeed necessary to specify a
selfadjoint extension for the Laplacian on a hyperbolic 3-manifold with Dehn
surgery type singularities. 
\begin{lemma}\label{not_ess_self}
Let $M$ be a hyperbolic 3-manifold with Dehn surgery type singularities. Then
$\Delta$ is not essentially selfadjoint. 
\end{lemma}
\begin{pf}
Since $\Delta_{max}$ (resp.~$\Delta_{min}$) is unitarily equivalent to
$\oplus_{\lambda \in \Lambda} (P_{V_{\lambda}})_{max}$ (resp.~to
$\oplus_{\lambda \in \Lambda} (P_{V_{\lambda}})_{min}$) on $T^2_{(0,R)}$,
it is enough to exhibit at least one $\lambda \in \Lambda$ such that
$P_{V_\lambda}$ is in the limit-circle case at $0$. 
Since $\lim_{r \rightarrow 0}r^2V_0(r)=-\frac{1}{4}$ we obtain by Theorem 2 in \cite{Sea1} that $P_{V_0}$ is indeed in the limit-circle case at $0$. 
Note that there might actually be infinitely many $\lambda \in \Lambda$ such
that $P_{V_\lambda}$ is in the limit-circle case at $0$. 
\end{pf}\\
\\
It is again a result of D.B.~Sears, cf.~\cite{Sea2}, that any self-adjoint extension of $P_V$ on $(0,R)$ has discrete spectrum if
$$
V(r) \geq -\frac{1}{4r^2}+A
$$
for some constant $A \in \R$. Note that under the same condition $P_V$ is semibounded by Hardy's inequality.

\begin{lemma}\label{discreteness}
Let $M$ be a hyperbolic 3-manifold with Dehn surgery type singularities.
Then $\spec \Delta_{Fr}$ is discrete. 
\end{lemma}
\begin{pf}
Since the essential spectrum is unaffected by removing a compact submanifold
with boundary, cf.~Proposition 1 in \cite{Bae}, it is sufficient to show that
$\spec \Delta_{Fr}$ is discrete on $T^2_{(0,R)}$. 

We use that $\Delta_{Fr}$ is unitarily equivalent to $\oplus_{\lambda \in
  \Lambda} (P_{V_{\lambda}})_{Fr}$ on $T^2_{(0,R)}$. From Theorem 1 in \cite {Sea2} we obtain that for $\lambda \in \Lambda$ any
  selfadjoint extension of $P_{V_\lambda}$ in $L^2(0,R)$ has discrete
  spectrum, hence in particular the Friedrichs extension. 
Observe that $P_{V_0}$ is a nonnegative operator, i.e.~$\langle
  P_{V_0}f,f\rangle_{L^2(0,R)} \geq 0$ for all $f \in C_0^{\infty}(0,R)$, and
  that $P_{V_\lambda}-P_{V_0} \geq C |\lambda|^2$ for all $\lambda \in
  \Lambda$ and some constant $C=C(R)>0$. 
We may estimate

\begin{align*}
\langle P_{V_\lambda} f,f \rangle_{L^2(0,R)} &= \langle P_{V_0} f,f
\rangle_{L^2(0,R)} + \langle P_{V_\lambda}-P_{V_0} f,f  \rangle_{L^2(0,R)}\\  
& \geq C|\lambda|^2 \|f\|^2_{L^2(0,R)}
\end{align*}
for $f \in C_0^{\infty}(0,R)$. For $l>0$ we obtain that $\spec
(P_{V_\lambda})_{Fr} \cap [0,l] \neq \emptyset$ only for finitely many
$\lambda \in \Lambda$, hence that $\spec \Delta_{Fr}$ is discrete on
$T^2_{(0,R)}$.  
\end{pf}\\

We finish this section with another application of the variational
characterization of eigenvalues.
\begin{lemma}\label{lem_bilip}
Let $X$ be a compact manifold with boundary, let $\beta>0$ and let $g_1$ and
$g_2$ be two 
Riemannian metrics on $X$ being $(1+\beta)$-quasi isometric:
\[ 
\frac{1}{(1+\beta)^2}\,g_2\le g_1 \le (1+\beta)^2\,g_2.
\]
For $j=1,2$ let $0<\lambda_1^j\le\ldots\lambda_i^j\le\ldots$ denote the eigenvalues of the Dirichlet problem on $(X,g_j)$, resp.~$0=\lambda_0^j<\lambda_1^j\le\ldots\lambda_i^j\le\ldots$ the eigenvalues of the Neumann problem.
Then for $i\ge 1$ one has
\[
\frac{1}{(1+\beta)^2}\,\lambda^1_i\le \lambda^2_i \le (1+\beta)^2\,\lambda^1_i.
\]
\end{lemma}
\begin{pf} 
Consider $f\in C^\infty(X)$. For the Rayleigh quotients one observes
\[
\frac{1}{(1+\beta)^2}\cdot \frac{\|df \|^2_{L^2(X,g_1)}}{\|f \|^2_{L^2(X,g_1)}}
\le \frac{\|df \|^2_{L^2(X,g_2)}}{\|f \|^2_{L^2(X,g_2)}} \le
(1+\beta)^2\cdot \frac{\|df \|^2_{L^2(X,g_1)}}{\|f \|^2_{L^2(X,g_1)}},
\]
and Theorem~\ref{varcharacterization} yields the claim for Neumann boundary conditions. Using the above estimate on the space of smooth functions compactly supported in the interior of $X$, then again Theorem~\ref{varcharacterization} yields the claim for Dirichlet boundary conditions.
\end{pf}

\section{Clustering}
Let $M^\infty$ be a hyperbolic 3-manifold of finite volume with a single
cusp in the following. For $M \in \Defo(M^\infty)$ we fix the
Friedrichs extension $\Delta_{Fr}$ for the Laplacian on $M$.

For a fixed interval $I \subset \R$ we view the spectral counting function
$\spc_{\Delta, M}I$ as a function on $\Defo(M^\infty)$. Note that for $\mu>0$ fixed the tube radius $R$ may also be considered as a function on $\Defo(M^\infty)$, cf.~Chapter \ref{Dehn_surgery}.  

\begin{thm}\label{spec_count}
The spectral counting function on $\Defo(M^\infty)$ satisfies 
$$
\spc_{\Delta,M}[1,1+x^2] = \frac{x}{\pi}R + O_x(1)
$$
for $x>0$. Here $O_x(1)$ denotes a function on $\Defo(M^\infty)$ which is
bounded in a neighbourhood of $M^\infty$. (This neighbourhood and the bounds
may depend on $x>0$.)
\end{thm}
We may reformulate this in terms of the tube shape, i.e.~the geometry of the
base of the singular tube:
\begin{cor}
The spectral counting function on $\Defo(M^\infty)$ satisfies 
$$
\spc_{\Delta,M}[1,1+x^2] = \frac{x}{2\pi}\left(\frac{1}{area\, T^2}\right) + O_x(1)
$$
for $x>0$.
\end{cor}
\begin{pf}
The assertion follows from Theorem \ref{spec_count} using Lemma \ref{tube_geom}.
\end{pf}\\
\\
We will bound the spectral counting function from above and from below,
i.e.~Theorem \ref{spec_count} will follow from Lemma \ref{lower_bound} and
Lemma \ref{upper_bound}:  

\begin{lemma}\label{upper_bound}
The spectral counting function on $\Defo(M^\infty)$ satisfies 
$$
\spc_{\Delta,M}[1,1+x^2]\leq\frac{x}{\pi}R + O_x(1)
$$
for $x>0$. Furthermore, $\spc_{\Delta,M}[0,1]=O(1)$.
\end{lemma}
\begin{pf}
Let $x \geq 0$. We claim that
$$
\spc_{\Delta,M}[0,1+x^2]\leq\frac{x}{\pi}R + O_x(1)\,,
$$
from which the assertions trivially follow. To prove this, we decompose 
$$
M = M_{[\mu,\infty)} \cup T^2_{(0,R]} \,,
$$
such that by the domain decomposition principle, cf.~\cite[Prop.~3]{Bae}, we get
$$
\spc^{Fr}_{\Delta,M}[0,1+x^2] \leq
\spc^{nat}_{\Delta,M_{[\mu,\infty)}}[0,1+x^2] +
  \spc^{nat}_{\Delta,T^2_{(0,R]}}[0,1+x^2]\,.   
$$

Here we choose {\em natural} boundary conditions for the Laplacian on the
various pieces, which are manifolds with boundary (compact or not),
cf.~\cite{Bae}. 
More precisely, for a manifold with boundary $N$ (compact or
not) we are considering the selfadjoint extension of $\Delta$ obtained by applying the
Friedrichs construction to the form domain $C_0^\infty(N)$, i.e.~compactly supported smooth functions on $N$, whose support is allowed to hit $\partial N$. If $N$ is in fact
compact, this amounts to choosing Neumann conditions on $\partial N$.  

We claim:
\begin{enumerate}
\item $\spc^{nat}_{\Delta,M_{[\mu,\infty)}}[0,1+x^2] = O_x(1)$.
\item $\spc^{nat}_{\Delta,T^2_{(0,R]}}[0,1+x^2] \leq \frac{x}{\pi}R + O_x(1)$.
\end{enumerate}

The first claim follows from the fact that the metrics on the thick parts are
uniformly quasi-isometric together with Lemma~\ref{lem_bilip}.

To prove the second claim, we will first show that only the $0$-mode contributes to
clustering. 
Towards that end, we choose $c=c(x)>0$ and a neighbourhood $W$ of
$M^\infty$ in $\Defo(M^\infty)$ such that $V_\lambda - V_0> 1+x^2$ on
the interval $(0,R-c]$ for all $M \in W$ and for all $\lambda \neq 0$. 
This is possible since there exists a constant $C>0$ independent of $\lambda
\neq 0$
such that for $r \in (0,R]$
$$
V_\lambda(r) - V_0(r) \geq Ce^{2(R-r)}(V_\lambda(R)-V_0(R))
$$
and hence for $r \in (0,R-c]$
$$
V_\lambda(r) - V_0(r) \geq Ce^{2c}(V_\lambda(R)-V_0(R))\,.
$$
Now $V_\lambda(R)-V_0(R)$ is bounded from below by the first positive
eigenvalue of $\partial T^2_{(0,R)}=\partial M_{[\mu,\infty)}$, which in
turn is bounded from below by some positive constant, since the induced metrics on
$\partial M_{[\mu,\infty)}$ are uniformly quasi-isometric in a neighbourhood of $M^{\infty}$ in $\Defo(M^{\infty})$. We choose $c>0$ accordingly and we may further decompose 
$$
T^2_{(0,R]} = T^2_{(0,R-c]} \cup T^2_{[R-c,R]}
$$
to obtain
$$
\spc^{nat}_{\Delta,T^2_{(0,R]}}[0,1+x^2] \leq
\spc^{nat}_{\Delta,T^2_{(0,R-c]}}[0,1+x^2] +
\spc^{nat}_{\Delta,T^2_{[R-c,R]}}[0,1+x^2]\,.
$$
Since the metrics on the regions $T^2_{[R-c,R]}$ are uniformly quasi-isometric in a neighbourhood of $M^{\infty}$ in $\Defo(M^{\infty})$, we find  
$$\spc^{nat}_{\Delta,T^2_{[R-c,R]}}[0,1+x^2]=O_x(1)\,.
$$ 
Clearly
$$
\spc^{nat}_{P_{V_\lambda},(0,R-c]}[0,1+x^2]=0
$$
since $V_\lambda -V_0 > 1+x^2$ on $(0,R-c]$ for $\lambda \neq 0$ and
$(P_{V_0})_{nat} \geq 0$, cf.~the proof
of Lemma \ref{discreteness}.

To estimate the contribution of the 0-mode, let $r_0>0$ be the unique positive zero of
$V_0(r)=2-\coth(2r)^2$. We may assume that $R-c>r_0$.  
Using the fact that $V_0 \geq 0$ for $r\geq r_0$ we obtain
\begin{align*}
\spc^{nat}_{P_{V_0},(0,R-c]}[0,1+x^2] &\leq
\spc^{nat}_{P_{V_0},(0,r_0]}[0,1+x^2] + \spc^{nat}_{P_{V_0},[r_0,R-c]}[0,1+x^2]\\ 
&\leq \spc^{Neu}_{P_{V_0},[r_0,R-c]}[0,1+x^2] +O_x(1)\\
&\leq \spc^{Neu}_{-\partial_r^2,[r_0,R-c]}[0,1+x^2] + O_x(1) = \dfrac{x}{\pi}R +
O_x(1)\,.
\end{align*}
This finishes the proof.
\end{pf}

\begin{lemma}\label{lower_bound}
The spectral counting function on $\Defo(M^\infty)$ satisfies 
$$
\spc_{\Delta,M}[1,1+x^2] \geq \frac{x}{\pi}R + O(1)
$$
for $x>0$.
\end{lemma}
\begin{pf}
Note that $V_0(r) \leq 1$ for $r >0$ and that $(P_{V_0})_{Fr} \geq 0$, cf.~the proof
of Lemma \ref{discreteness}. By domain monotonicity, cf.~\cite[Prop.~2]{Bae},
and the variational characterization of eigenvalues we get 
\begin{align*}
\spc^{Fr}_{\Delta,M}[0,1+x^2] &\geq \spc^{Fr}_{\Delta,T^2_{(0,R)}} [0,1+x^2]\\
& \geq \spc^{Fr}_{P_{V_0},(0,R)}[0,1+x^2]\\
& \geq \spc^{Fr}_{-\partial_r^2+1,(0,R)}[0,1+x^2]\\
& = \spc^{Dir}_{-\partial_r^2,(0,R)}[0,x^2] = \frac{x}{\pi} R +O(1)\,.
\end{align*}
Now by Lemma \ref{upper_bound} we
have $\spc_{\Delta,M}[0,1] = O(1)$, hence we obtain, as claimed, that $\spc_{\Delta,M}[1,1+x^2]
\geq \frac{x}{\pi}R + O(1)$.
\end{pf}

\section{Convergence of the small eigenvalues}
 
Let $M^\infty$ be a hyperbolic 3-manifold of finite volume with a single
cusp in the following. 
Let $0=\lambda^\infty_0< \lambda_1^{\infty} \leq \ldots \leq
\lambda^\infty_{k^\infty} < 1$ denote the eigenvalues of the Laplacian on
$M^\infty$ below the essential spectrum. 
We fix some $\Lambda<0$ with $\lambda_{k^\infty}<\Lambda< 1$.
For $M \in \Defo(M^\infty)$ let $0=\lambda_0 < \lambda_1 \le \ldots\le\lambda_k\le
\ldots\to\infty $ denote the eigenvalues of $\Delta_{Fr}$ on $M$, and let 
$k(\Lambda)$ denote the number of positive eigenvalues strictly below
$\Lambda$, which means $\lambda_{ k(\Lambda)}<\Lambda\le \lambda_{ k(\Lambda)+1}$. 
\begin{thm}\label{convergence}
For each $\varepsilon > 0$ there exists a neighbourhood $W$ of $M^\infty$ in $\Defo(M^\infty)$ such that for $M \in W$ one has: 
\begin{enumerate}
\item $k(\Lambda)= k^\infty$,
\item $|\lambda_i - \lambda_i^\infty| < \varepsilon$ for $i=1, \ldots,k^\infty$.
\end{enumerate}
\end{thm}
We may reformulate this in terms of generalized Dehn surgery coefficients:
\begin{cor}
For each $\varepsilon > 0$ there exists a neighbourhood $W'$ of $\infty$ in $(\R^2 \cup \{ \infty \})/\pm 1$  such that for $(x,y)\in W'$ one has for
$M=M^{(x,y)}$:   
\begin{enumerate}
\item $k(\Lambda)= k^\infty$,
\item $|\lambda_i - \lambda_i^\infty| < \varepsilon$ for $i=1,
  \ldots,k^\infty$.
\end{enumerate}
\end{cor}
\begin{pf}
The assertion follows from Theorem \ref{convergence} together with Theorem
\ref{dehn_surgery}. 
\end{pf}\\
\\
In order to prove Theorem~\ref{convergence} we will extend Colbois' and
Courtois' method of {\em parties fuyantes} (see \cite{CC}) to hyperbolic
manifolds with Dehn surgery type singularities.
We need some preliminary considerations.

If $\mu>0$ is small enough the $\mu$-thin part of $M^\infty$ consists of a
rank-2 cusp: $M^\infty_{(0,\mu)}$ can isometrically be identified with
$(0,\infty)\times T^2$ carrying the Riemannian metric $dt^2+e^{-2t}g_{T^2}$,
and $M^\infty_{[\mu,\infty)}$ is a smooth compact manifold with boundary.
Let $0<\lambda_0^\infty(\mu)<\lambda_1^\infty(\mu)
\le\ldots\le\lambda_i^\infty(\mu)\le\ldots\to\infty$
be the eigenvalues of the Dirichlet problem on $M^\infty_{[\mu,\infty)}$,
and denote the index of the last eigenvalue strictly below $1$ by
$k^\infty_\mu$.  
Then the same proof as the one of (1.4) in \cite{CC} shows:

\begin{lemma}\label{lem:ausschoepf}
For each $\varepsilon>0$ with $\lambda^\infty_{k^\infty}+\varepsilon < 1$ there is a
$\mu_\varepsilon>0$ such that for all $\mu<\mu_\varepsilon$ one has:
\begin{enumerate}
\item $k^\infty_\mu = k^\infty$,
\item $\lambda_i^\infty\le\lambda_i^\infty(\mu)<\lambda_i^\infty+\varepsilon$ for $i=0,
  \ldots,k^\infty$. 
\end{enumerate}
And, by domain monotonicity, one also has $\lambda_{k^\infty+1}^\infty(\mu)\ge 1$.
\hfill$\boxbox$
\end{lemma}

Next we note that elliptic regularity is a local statement, hence it also
holds on manifolds with Dehn surgery type singularities.  
We get that any eigenfunction~$f\in L^2(M)$ of $\Delta_{Fr}$ is smooth.

For the moment let us fix $\rho>0$ and $c>4$ and let us consider a hyperbolic
3-manifold $M$ containing a singular tube $T^2_{(0,\rho+c)}$ 
isometric to $(0,\rho+c)\times T^2$ with the metric given in
(\ref{singtubemetric}).  
For $0\le a<b\le \rho+c$ we denote the domain that
corresponds to $(a,b]\times T^2$ 
by  $T^2_{(a,b]}$ and for $0<r\le \rho+c$ we denote the hypersurface that corresponds to $\{r\}\times T^2$ by $T_r^2$.
\begin{lemma}\label{rlemma}
For any $\Lambda>0$ and for any $f\in C^\infty(M)$ with
$\|f\|^2_{L^2(M)}=1$ and $\|df\|^2_{L^2(M)}<\Lambda$ there exists an $r\in
[\rho+2,\rho+c]$ 
such that 
\begin{enumerate}
\item $\|f\|^2_{H^1(T^2_{(r-1,r]})}<\frac{2\cdot(1+\Lambda)}{c-4}$ and
\item $\int_{T_r^2}\left\{|f|^2+|df|^2 \right\}<\frac{2\cdot(1+\Lambda)}{c-4}$.
\end{enumerate}
\end{lemma}
\begin{pf}
In fact the proof is the same as the one of \cite[Lemme~2.4]{CC}:
One defines $F=|f|^2+|df|^2$ and gets
\[
1+\Lambda>\int_M F\ge \int_{T^2_{(\rho+2,\rho+c]}}F
\ge\sum_{k=1}^{[\frac{c}{2}]-1}\int_{T^2_{(\rho+2k,\rho+2k+2]}}F. 
\]
Hence there is a $k\in\left\{1,\ldots,[\frac{c}{2}]-1\right\}$ with
\[
\frac{1+\Lambda}{[\frac{c}{2}]-1}>\int_{T^2_{(\rho+2k,\rho+2k+2]}}F
  =\int_{\rho+2k}^{\rho+2k+2}dr\int_{T_r^2}F.   
\]
Furthermore there is an $r\in[\rho+2k+1,\rho+2k+2]\subset [\rho+2,\rho+c]$
with $\int_{T_r^2}F<\frac{1+\Lambda}{[\frac{c}{2}]-1}
\le\frac{2\cdot(1+\Lambda)}{c-4}$ and we are done.
\end{pf}
\begin{cor}
For $\Lambda>0$, $f\in C^\infty(M)$ and $r\in[\rho+2,\rho+c]$ as in
Lemma~\ref{rlemma} with $\|\Delta f\|_{L^2}<\infty$
one has 
\begin{equation}\label{eq:green}
\left|\int_{T^2_{(0,r]}}|df|^2-\int_{T^2_{(0,r]}}f\cdot\Delta f \right|
    < \frac{2\cdot(1+\Lambda)}{c-4}.
\end{equation}
\end{cor}
\begin{pf}
We use the Green's formula in Corollary~\ref{L2_green} to get
\[
\left|\int_{T^2_{(0,r]}}|df|^2-\int_{T^2_{(0,r]}}f\cdot\Delta f \right| =
\left|\int_{T_r^2}f\cdot\partial_r f \right|\le
\sqrt{\int_{T_r^2} f^2}\cdot
\sqrt{\int_{T_r^2}(\partial_r f)^2}.
\]
Then Lemma~\ref{rlemma} finishes the proof.
\end{pf}

\begin{lemma}\label{est:lower_tube}
Let $\Omega=T^2_{(0,R)}$ be a singular tube of radius $R$. 
Consider the Friedrichs extension $\Delta_{Fr}$ of the Laplace operator acting
on $C^\infty_0\left(\Omega \right)$.
For $h\in\dom(\Delta_{Fr}) $ one has
\[\|h\|_{L^2(\Omega)} \le \|dh\|_{L^2(\Omega)}.\]
\end{lemma}
\begin{pf}
For $h\in C^\infty_0\left(\Omega \right)$ one gets
\begin{equation}\label{CSU}
\int_{\Omega}\left|\partial_r(h^2) \right|
=\int_{\Omega}\left|2h\cdot\partial_rh \right|
\le 2\| h\|_{L^2(\Omega)}\cdot \| dh\|_{L^2(\Omega)}.
\end{equation}
As $h$ is compactly supported we may integrate by parts and obtain
\begin{align*}
\| h\|_{L^2(\Omega)}^2
&=\int_{T^2}d\theta\, dz\int_0^R h^2 \cdot
\frac{d}{dr}\left(\frac{1}{2}\sinh^2(r) \right)dr \\ &
=-\int_{T^2}d\theta\, dz\int_0^R 
\frac{1}{2}\sinh^2(r)  \cdot \partial_r\left(h^2\right)dr\\
&=-\tfrac{1}{2}\int_\Omega\tanh(r) \partial_r\left(h^2\right) 
\le\tfrac{1}{2}\int_\Omega \left|\partial_r(h^2) \right|
\le \| h\|_{L^2(\Omega)}\cdot \| dh\|_{L^2(\Omega)}
\end{align*}
by (\ref{CSU}).
For the Rayleigh quotient this means 
\[ \frac{\| dh\|_{L^2(\Omega)}^2}{\| h\|_{L^2(\Omega)}^2}\ge 1.\]
By the construction of Friedrichs extensions this lower bound holds for
any $h\in\dom(\Delta_{Fr})$, and we are done.
\end{pf}
\begin{cor}
Let $M$ be a hyperbolic 3-manifold containing a singular tube
  $T^2_{(0,\rho+c]}$, where  $\rho>0$ and $c>4$. 
Then for any $f\in\dom(\Delta_{Fr})$ and any $r\in[\rho+2,\rho+c]$ one gets 
\begin{equation}\label{green_error}
\int_{T^2_{(0,r]}}|f|^2 -\int_{T^2_{(0,r]}} |df|^2 \le 3\cdot
    \|f\|^2_{L^2(T^2_{(r-1,r]})}. 
\end{equation}
\end{cor}
\begin{pf}
We choose a cut-off function $u\in C^\infty(M)$ with $\supp(u)\subset T^2_{(0,r)}$,
$0\le |u|\le 1$,  $|\grad(u)|\le 2$ and $u|_{T^2_{(0,r-1]}}\equiv 1$
and consider $h=u\cdot f$.
Then $h$ is in the domain of the Friedrichs extension of the Laplacian acting
on $C^\infty_0(T^2_{(0,r)})$, and by Lemma~\ref{est:lower_tube} we get
\begin{align*}
\int_{T^2_{(0,r]}}|f|^2-\|f\|^2_{L^2(T^2_{(r-1,r)})}& \le
\int_{T^2_{(0,r]}}|h|^2\le
\int_{T^2_{(0,r]}}|dh|^2\\
&\le
\int_{T^2_{(0,r]}}|df|^2+2\cdot\|f\|^2_{L^2({T^2_{(r-1,r]}})},
\end{align*}
which finishes the proof.
\end{pf}\medskip

Now let $f\in C^\infty(M)$ be a normalized eigenfunction of $\Delta_{Fr}$ corresponding
to an eigenvalue $\lambda<\Lambda<1$, which means $\Delta_{Fr}f=\lambda f$ and
$\|f\|_{L^2(M)}=1$.
We take $r\in[\rho+2,\rho+c]$ as in Lemma~\ref{rlemma}, and
we abbreviate 
$$
A:=\int_{T^2_{(0,r]}}|df|^2 \,,\quad B:=\int_{T^2_{(0,r]}} |f|^2 \,,\quad \eta:=\frac{2\cdot(1+\Lambda)}{c-4}.
$$
One gets $|A-\lambda\cdot B|<\eta$ by (\ref{eq:green}) and $B-A < 3\eta$ by (\ref{green_error}) and Assertion 1 of Lemma \ref{rlemma}. Together this yields $B-\lambda\cdot B<A+3\eta-\lambda B<4\eta$ and hence
$B<\frac{4\eta}{1-\lambda}<\frac{4\eta}{1-\Lambda}$. 
Using $A+B \leq (1+\Lambda)B + \eta$ we conclude:

\begin{prop}\label{prop:sobotube}
Let~$f\in C^\infty(M)$ be an eigenfunction of~$\Delta_{Fr}$ corresponding to
an eigenvalue~$\lambda<\Lambda<1$ with~$\|f\|_{L^2}=1$.
Then there is an $r\in[\rho+2,\rho+c]$ such that
\[
\|f \|_{H^1(T^2_{(0,r]})}^2<\frac{40}{1-\Lambda}\cdot\frac{1}{c-4}. 
\]  
\end{prop}
Finally we have all ingredients to prove Theorem \ref{convergence}:\\

\noindent
{\textit{Proof of Theorem~\ref{convergence}.}}
Let $\varepsilon > 0$. Without loss of generality we may assume that
$\lambda^\infty_{k^\infty}+\varepsilon <\Lambda$ and $\Lambda+\varepsilon<1$. Applying Lemma~\ref{lem:ausschoepf} for $\varepsilon/2$ we obtain $\mu_{\frac{\varepsilon}{2}}>0$ such
that for any $\mu<\mu_{\frac{\varepsilon}{2}}$ one has
$k^\infty_\mu=k^\infty$ and
\begin{align}
\lambda^\infty_i&\le \lambda^\infty_i(\mu)\le
\lambda^\infty_i+\frac{\varepsilon}{2}\qquad\mbox{ for
}i=0,\cdots,k^\infty,\mbox{ and}\label{est:inftymu}\\
 1&\le \lambda^\infty_{k^\infty+1}(\mu).\label{est:inftymu2}
\end{align}
We choose $\mu\in(0,\mu_{\frac{\varepsilon}{2}})$. 

For $M \in \Defo(M_\infty)$ we denote the
Dirichlet-eigenvalues of $M_{[\mu,\infty)}$ by $\lambda_i(\mu)$.
By Lemma \ref{lem_bilip} there exists $\beta>0$ such that if the $\mu$-thick parts $M_{[\mu,\infty)}$ and
  $M^\infty_{[\mu,\infty)}$ are $(1+\beta)$-quasi-isometric, then one has
\begin{align}
\lambda^{\infty}_i(\mu)-\frac{\varepsilon}{2} &\le
\lambda_i(\mu)\le
\lambda^{\infty}_i(\mu)+\frac{\varepsilon}{2} \qquad\mbox{ for
}i=0,\cdots,k^\infty,\mbox{ and}\label{contradiction:eps} \\
\lambda_{k^\infty+1}(\mu) &>
\Lambda+\frac{\varepsilon}{2}.\label{contradiction:eps2} 
\end{align}
From (\ref{est:inftymu}) and
(\ref{contradiction:eps}) we get for $i=0,\ldots,k^\infty$ 
\begin{equation}\label{est:nearlyfinished}
\lambda^\infty_i-\frac{\varepsilon}{2}\le
\lambda_i(\mu)\le
\lambda^\infty_i+\varepsilon.
\end{equation}
Domain monotonicity applied to
$M_{[\mu,\infty)} \subset M$ gives 
\begin{equation}\label{est:nearlynearly}
 \lambda_i\le\lambda_i(\mu)\le \lambda^\infty_i+\varepsilon<\Lambda
\end{equation}
for $i=0,\ldots,k^\infty$, and therefore $k(\Lambda)\ge k^\infty$.

For $\lambda_0,\ldots,\lambda_{k^\infty}$ we take $L^2$-orthonormal
eigenfunctions $f_0,\ldots,f_{k^\infty}$.
Now for $i\le k^\infty$ Proposition~\ref{prop:sobotube} applies and we find
$r_i\in[\rho+2,\rho+c)$ with 
\[ 
\|f_i\|^2_{H^1(T^2_{(0,r_i]})}<\frac{40}{1-\Lambda}\cdot\frac{1}{c-4}.
\]
We choose $\varphi_i\in C^\infty(M)$ with $0\le\varphi_i\le 1$,
$\varphi_i\big|_{T^2_{(0,\rho]}}\equiv 0$, $\varphi_i\big|_{M\setminus
  T^2_{(0,r_i]}}\equiv 1$ and $|\grad \varphi_i|\le 2$ and get
\begin{equation}\label{H1close}
\|f_i-\varphi_i\cdot f_i\|^2_{H^1(M)}
\le 3\cdot \|f_i\|^2_{H^1(T^2_{(0,r_i]})}<\frac{120}{1-\Lambda}\cdot\frac{1}{c-4}.
\end{equation}
Now $g_i=\varphi_i\cdot f_i$ are smooth functions with compact support inside the interior of $M_{[\mu,\infty)}$. By (\ref{H1close}) we get that for $c$ large enough the functions $g_i$ are linearly independent and their Rayleigh quotients are arbitrarily close to those of the $f_i$ (the difference being uniformly controlled by $c$). Hence, by choosing $c$ larger than some constant only depending on the number $k_{\infty}$, we can achieve $\lambda_i(\mu)\le
\lambda_i+\frac{\varepsilon}{2}$ for $i=0,\ldots,k^\infty$. Together with (\ref{est:nearlyfinished}) and (\ref{est:nearlynearly}) this means
\[
\lambda^\infty_i-\varepsilon\le
\lambda_i\le \lambda^\infty_i+\varepsilon\qquad\mbox{ for }i=0,\ldots,k^\infty .
\]
If we suppose $\lambda_{k^\infty+1}<\Lambda$, the same argument as above (possibly involving a new choice of $c$) yields 
$\lambda_{k^\infty+1}(\mu)<\Lambda +\frac{\varepsilon}{2}$ which is
a contradiction to (\ref{contradiction:eps2}).
Therefore $k(\Lambda)\le k^\infty$.

For these chosen $\mu,\beta,c$  we take the neighbourhood $W$ of
$M^\infty\in\Defo(M^\infty)$ which is given in 
Lemma~\ref{lem_minftydef}, and by the preceding arguments it is clear that $W$ is the desired
neighbourhood in Theorem~\ref{convergence}.\hfill$\boxbox$


\end{document}